\documentstyle[11pt]{article}
\newtheorem{thm}{Theorem}

\newtheorem{conj}[thm]{Conjecture}

\newtheorem{question}[thm]{Question}

\newcommand{\beq}[1]{\begin{equation}\label{#1}}
\newcommand{\enq}[0]{\end{equation}}

\newcommand{\qed}[0]{\begin{flushright} \rule{2mm}{3mm} \end{flushright}}
\newcommand{\C}[2]{{{#1}\choose{{#2}}}}
\newcommand{\ga}[0]{\alpha }
\newcommand{\gb}[0]{\beta }

\newcommand{\gD}[0]{\Delta }

\newcommand{\eps}[0]{\varepsilon }

\newcommand{\0}[0]{\emptyset}

\newcommand{\A}[0]{{\cal A}}
\newcommand{\B}[0]{{\cal B}}

\newcommand{\f}[0]{{\cal F}}
\newcommand{\g}[0]{{\cal G}}
\newcommand{\h}[0]{{\cal H}}

\newcommand{\bn}[0]{\bigskip\noindent}
\newcommand{\mn}[0]{\medskip\noindent}
\newcommand{\nin}[0]{\noindent}

\newcommand{\sub}[0]{\subseteq}
\newcommand{\sm}[0]{\setminus}
\renewcommand{\dots}[0]{,\ldots,}

\newcommand{\ra}[0]{\rightarrow}
\newcommand{\Ra}[0]{\Rightarrow}

\newcommand{\E}[0]{{\sf E}}
\newcommand{\pe}[0]{p_{\E}}

\begin{document}
\renewcommand{\thefootnote}{\fnsymbol{footnote}}

\begin{center}
{\Large\bf Thresholds and expectation thresholds}
\end{center}

\begin{center}
J. Kahn\footnote{Supported by NSF grant DMS0200856.}
and G. Kalai\footnote{Supported by NSF and ISF grants.}\\
\end{center}
\date{}

\mn
{\em Monotone properties and thresholds}

\medskip
We use $2^X$ for the collection of subsets of a finite set $X$.
Given $X$ of size $n$ and $p\in [0,1]$, $\mu_p$ is the product
measure on $2^X$ given by $\mu_p(S) = p^{|S|}(1-p)^{n-|S|}$.
Recall that $\f\sub 2^X$ is often called a {\em property} of
subsets of $X$, and is said to be {\em monotone (increasing)} if
$B\supseteq A\in\f\Ra B\in\f$.
For a monotone $\f \sub 2^X$, denote by $p_c(\f)$ that $p\in
[0,1]$ for which $\mu_p(\f)$ ($=\sum\{\mu_p(S): S \in \f\}$)
$=1/2$. We will assume throughout this note that $\f$ is a
monotone property, and will always exclude the trivial cases
$\f=2^X$ and $\f=\0$, so in particular $p_c(\f)$ exists and is
unique.

Recall that for monotone $\f_n\sub 2^{X_n}$ ($n=1,2,\ldots$), a
function $p_0(n)$ is a {\em threshold} for the sequence $\{\f_n\}$
if
\begin {equation}
\mu_{p(n)}(\f_n)\ra\left\{\begin{array}{ll}
0&\mbox{if $p(n)/p_0(n)\ra 0$}\\
1&\mbox{if $p(n)/p_0(n)\ra \infty$.}
\end{array}\right.
\end {equation}
(See {\em e.g.}
\cite{AS}, \cite{BRG} or \cite{JLR}.
Of course $p_0(n)$ is not quite unique, but
following common practice we will often say
``the" threshold when we should really say ``a.")
It follows from \cite{BT}
that every $\{\f_n\}$ has a threshold, and that in fact
(see {\em e.g.}
\cite[Proposition 1.23 and Theorem 1.24]{JLR})
$p_0(n):=p_c(\f_n)$ is
a threshold for $\{\f_n\}$.
So the quantity $p_c$ conveniently captures threshold behavior,
in particular allowing us to dispense with sequences.

\begin{conj}\label{C1}
For any (nontrivial monotone) $\f \sub 2^X$ there exist $\g\sub
2^X$ and $q\in [0,1]$ such that

\mn {\rm (a)} each $B\in \f$ contains some member of $\g$,

\mn {\rm (b)} $\sum_{A\in \g}q^{|A|} < 1/2$, and

\mn
{\rm (c)}  $q>p_c(\f)/(K\log n)$,

\mn where $K$ is a universal constant and $n=|X|$.
\end{conj}
Note that, for a given $q$, existence of a $\g$ for which (a) and (b)
hold trivially gives $\mu_q(\f)<1/2$;
so the conjecture says that for any $\f$ as above
there is a trivial lower bound on $p_c(\f)$ which
is within a factor $O(\log n)$ of the truth.
To put this another way,
define $q(\f)$ to be the supremum
of those $q$'s for which there exists $\g$ satisfying (a) and (b).
Then $p_c(\f)\geq q(\f)$ is trivial, and Conjecture~\ref{C1}
would say $p_c(\f) <Kq(\f)\log n$.
It is easy to see (familiar examples will be recalled below)
that the factor $\log n$ cannot be improved.

(Condition (b) is also considered by Talagrand in
\cite{oldtal} and \cite{Talagrand}.
In particular a hope of \cite{oldtal} (which is also primarily a problem
paper) is to say that in certain geometrically-motivated situations
the threshold essentially coincides with its trivial lower bound
(more precisely, the gap is $O(1)$).
This is different from what we are asking, but it may be that
the underlying difficulties are similar.)


In the rest of this note
we will briefly mention a few consequences of,
and questions
related to, Conjecture~\ref{C1}.
As some of these indicate, the conjecture
is extremely strong.
It would probably be more
sensible to conjecture that it is {\em not} true;
but it does not seem easy to disprove, and
we think a counterexample would also be quite interesting.


\mn
{\em Graph properties}

\medskip
A basic result for random graphs states that a threshold for
the (usual) random graph $G_{n,p}$
to contain a given (fixed) subgraph $H$
is $n^{-1/m(H)}$, where
$m(H)=\max\{|E(H')|/|V(H')|:H'\sub H\}.$
(This was proved for ``balanced" graphs in \cite{ER}, and
for general graphs in \cite{Boll81}.)
Equivalently,
one may take as a threshold
the least $p=p(n)$ such that
for each $H'\sub H$ the expected number of copies of
$H'$ in $G_{n,p}$ is at least $1$.

This last $p(n)$---call it the {\em ``expectation threshold"}---still
makes sense if,
instead of a fixed $H$, we
consider a {\em sequence} $\{H_n\}$,
where, allowing isolated vertices,
we assume $|V(H_n)|=n$.
Formally, for an arbitrary $H$,
define $\pe(H)$ to be the least $p$ such that,
for every spanning $H'\sub H$,
$
(|V(H')|!/|Aut(H')|)p^{|E(H')|} \geq 1.
$
Then the $p(n)$ of the preceding paragraph is the same as
$\pe(H_n)$, where $H_n$ is gotten from $H$ by
adding $n-|V(H)|$ isolated vertices.

Of course in this more general situation, $\pe$ may no longer
capture the true threshold behavior. Well-known examples are
perfect matchings and Hamiltonian cycles in $G_{n,p}$, for each of
which the expectation threshold is (easily seen to be) on the
order of $1/n$, while the actual threshold is $\log n/n$. (This
was proved for matchings in \cite{ER}, and for Hamiltonian cycles
in \cite{Posa} and \cite{Korsh}. It's easy to see that $\log n/n$
is a lower bound, since existence of a perfect matching and
Hamiltonicity require minimum degree at least 1 and 2
respectively. The beautiful ``stopping time" versions of these
results (\cite{ER} for matchings, \cite{KS} and \cite{B} for
Hamiltonian cycles) say, in a precise way, that in each case
failure to satisfy the degree requirements is the main source of
failure to have the desired property.)

Now, for an $n$-vertex $H$,
write $p_c(H)$ and $q(H)$ for $p_c(\f)$ and $q(\f)$,
where $\f\sub 2^{\C{[n]}{2}}$ is the collection of graphs on $[n]$
that contain copies of $H$.
So again, $p_c(H_n)$ is a threshold function for ``$G_{n,p}\supseteq H_n$."
Moreover, for any $H$, $q(H)\geq \pe(H)/2$
(the 1/2 should of course be ignored).
To see this, notice that we have an alternate definition of
$\pe(H)$ analogous to that of $q(H)$:
it is the supremum of those $p$'s for which there exists
$H'\sub H$ so that, with $\g$ the collection of copies of $H'$,
$\sum\{p^{|A|}:A\in \g\}<1$.
So the following conjecture, which was actually our starting point,
contains Conjecture~\ref{C1} for $\f$'s of this (subgraph
containment) type.

\begin{conj}\label{C2}
For any $H$, $p_c(H) <K\pe(H)\log |V(H)|$,
where $K$ is a universal constant.
\end{conj}
On the other hand, it could be that Conjecture~\ref{C1} implies
Conjecture~\ref{C2}, as would follow from a positive answer to

\begin{question}
Is it true that, for some fixed K, $q(H)<K\pe(H)$ for every $H$?
\end{question}
In other words, is it true that in the present situation the value
of $q$ does not change significantly if we require $\g$ (in the
definition of $q$) to consist of all copies of a single graph?
Pushing our luck, we could extend this to a general graph property
$\f$ (meaning, as usual, that membership in $\f$ depends only on
isomorphism type):  can we in this case require that $\g$ also be
a graph property (without significantly affecting $q$)? In fact,
we don't even know that the answer to the following is negative.

\begin{question}
For $\f\sub 2^X$ let $q^*(\f)$ be the supremum of those $q$'s
for which there is some Aut($\f$)-invariant $\g$ satisfying
(a) and (b) of Conjecture~\ref{C1}.
Is it true that $q(\f)\leq Kq^*(\f)$ for some universal K?
\end{question}

Returning briefly to subgraph containment, it is not hard to see
that for an $n$-vertex tree $T$ with maximum degree $\gD(T)=\gD$,
both $\pe(T) $ and $q(T)$ are always between $\Omega(1/n)$ and
$O(\Delta/n)$, and are $\Theta(\Delta/n)$ if $\Delta=\Omega(\log
n)$. In particular, either of the above conjectures would imply
that for any fixed $\Delta$ we have $p_c(T)<O(\log n/n)$ for any
$T$ of maximum degree at most $\gD$. For unbounded degrees, the
conjectures are not as precise, but as an aside we may mention the
natural guess that there are constants $K_1>0$ and $K_2$ such that
for any tree $T$ (with $n$ vertices)
\begin {equation}
K_1\max\{\log n/n,\gD(T)/n\} < p_c(T) < K_2\max\{\log n/n,\gD(T)/n\}.
\end {equation}

Finally, we don't even know a counterexample to

\begin{conj}\label{C2prime}
Given $\eps >0$ there is a $K$ such that any $H$ with
$\eps <p_c(H)< 1-\eps$ satisfies
$p_c(H) <K\pe(H)$.
\end{conj}

\nin
{\em Hypergraph matching}

\medskip
Write $\h_k(n,p)$ for the random $k$-uniform hypergraph on vertex
set $[n]$ in which each $k$-set is an edge with probability $p$,
independent of other choices. A celebrated question first raised
and studied by
Schmidt and Shamir \cite {SS,Erdos} asks: for fixed $k$ and $n$
ranging over multiples of $k$, what is the threshold for
$\h_k(n,p)$ to contain a perfect matching (meaning, of course, a
collection of edges partitioning the vertex set)?
The best published progress to date, building in part on
\cite{FJ} and some antecedents, is due to
J.H. Kim \cite{Kim}, who shows that the threshold is at most
$O(n^{-\ga_k})$ with
$\ga_k = k-1-1/(5+2(k-1))$;
and very recently A. Johansson \cite{Joh} announced that he can
improve this to $O(n^{-k+1+o(1)})$.
But it is natural to expect
that, as for
$G_{n,p}$, one usually has a perfect matching as soon as there are
no isolated vertices, which
would say in particular that
(this was perhaps first proposed explicitly in \cite{CFMR})
the true threshold is $n^{-k+1}\log n$.
In fact this would follow  from Conjecture~\ref{C1},
since it's easy to see that for
$\f$ the collection of $k$-uniform hypergraphs on $[n]$
containing perfect matchings
(a property of subsets of $X=\C{[n]}{k}$),
one has $q(\f) =\Theta(n^{-k+1})$.

A 
related problem asks for the
threshold for $G_{n,p}$ to contain a ``triangle factor"
(collection of triangles partitioning the vertex set).
Here again one expects that isolated vertices (i.e.
vertices not contained in triangles) are the main
obstruction, which would make the threshold
$n^{-2/3}\log^{1/3}n$.  In this case Conjecture~\ref{C1} is
a little off, saying only that the threshold is at most
$n^{-2/3}\log n$; but this is still much better than what's presently
known. (Here again the best published bound is from \cite{Kim}
($n^{-2/3+1/18}$, improving on
\cite{Kriv}), and
$n^{-2/3 +o(1)}$ has been announced in \cite{Joh}.)

\mn
{\em Isoperimetry}

\medskip
We would like to suggest a possible, if speculative,
connection between the preceding questions and isoperimetric
behavior.
The most obvious difficulty posed by Conjecture~\ref{C1} is that
of identifying the requisite $\g$.  While such an identification based
purely on combinatorial information seems difficult,
a beautiful result of Friedgut \cite{Fri} does achieve
something of the sort using the Fourier transform
(see the comments following Conjecture~\ref{C3}).
This somewhat motivates the approach proposed here.

Given a monotone $\f \subset 2^X$, let $m(p) =
\mu_p(\f).$ The derivative $m'(p)$
is, according to a fundamental
observation of Russo (\cite{Russo} or e.g. \cite {Grimmett}), equal to the ``total
influence'' of $\f$ w.r.t. $\mu_p$, namely
$$
I=I_p(\f)=\sum_{i=1}^n \mu_p(|\f\cap \{S, S\triangle \{i\}\}|=1).
$$
This measure of the edge boundary of $\f$ has recently been of
importance in a number of contexts; see e.g.
\cite{KKL},\cite{Fri},\cite{KaS}. (The combination of Russo's
lemma and lower bounds on influence has been a powerful tool for
proving sharp threshold behavior. A few recent results have
also used such sharp threshold information as a tool in estimating
{\em location} of thresholds; see e.g. \cite{ANP},\cite{BR}.)

In the present context the edge isoperimetric inequality takes
the form

\begin {equation}\label{iso}
p \cdot I_p(\f) \ge    m(p) \log _p  m(p).
\end {equation}
Equality holds whenever $\f$ is a subcube, $\{S \subset X: S
\supseteq R\}$ for some $R$. Though one would think (\ref{iso})
would be known, we could not find it in the literature, so
indicate a proof at the end of this note. A derivation, based on
log-Sobolev inequalities, of something containing
a slightly weaker version of (\ref{iso}) (in which
the right hand side
is multiplied by a constant
$C_p$) is given in \cite[p.111, (5.28)]{Ledoux}.

For $\g\sub 2^X$, write $\langle \g\rangle$ for the monotone
property generated by $\g$; that is, $\langle \g\rangle = \{B\sub
X: \exists A\in \g, A\sub B\}$. (So, for example, (a) of
Conjecture~\ref{C1} is ``$\f\sub\langle \g\rangle$.") We would
like to say, very roughly, that if (\ref{iso}) is fairly tight for
a pair $(\f,p)$, then $\f$ is close (in $\mu_p$) to
some property
generated by small sets.
We will say that $\f$ is $(C,p)$-{\em optimal} if (\ref{iso})
is tight to within the factor $C$; that is,
$$
p \cdot m ' (p) \le C m(p)  \log_p m(p).
$$
Thus subcubes are $(1,p)$-optimal for every $p$.

\begin {conj}
\label {C3}
Given $C$, there are $K,\delta >0$ such that
for any $(C  \log (1/p) ,p)$-optimal
$\f$ (and $m(p)=\mu_p(\f)$):

\mn
{\rm (a)} There is an $R \subset X$ of size at most
$K\log (1/ m(p)), $ such that
$$
\mu_p (S \in \f | S\supseteq R) \ge (1+\delta) m(p).
$$

\mn
{\rm (b)} There is a family $\g \subset \{A\sub X: |A|\leq
K  \log (1/ m(p))\}$ such that
$\mu_p(\f\triangle\langle\g\rangle)< o(m(p))$.

\mn
Furthermore, there is some fixed $C>0$ such that for any
$(C  \log (1/p) ,p)$-optimal $\f$,

\mn
{\rm (c)} there is $\g \subset 2^X$ such that
$\f\sub\langle\g\rangle$ and $\sum \{p^{|S|}:S \in {\g}\} \le
1/2.$
\end {conj}

For $\mu_p(\f)$ bounded away from 0 and 1, (b) is a conjecture of
Friedgut \cite {Fri}, which he proved in case $\f$ is a graph
property, and (a) is similar to a theorem of Bourgain \cite{Bou}.
%
The extra $\log(1/p)$ in Conjecture~\ref{C3} seems a little
unnatural, but makes sense in light of whatever examples we've looked at
(see in particular the prototypical ``dual tribes" example below).
There is also some vague feeling that in ``cases of interest"
one can multiply the right hand side of (\ref{iso}) by $\Omega(\log(1/p))$.
(Of course Conjecture~\ref{C3} may be regarded as a statement to
this effect.
Also, for instance, the extra log is not needed if
 $\mu_p(\f)$ is bounded
away from 0 and 1, since (\ref{iso}) applied
to the {\em dual} family,
$\f^*= \{A\sub X: X\sm A\not\in \f\}$, gives (in general)
$$
(1-p)I\geq (1-\mu_p(\f))\log_{1-p}(1-\mu_p(\f))
$$
(note $\mu_{1-p}(\f^*) = 1-\mu_p(\f)$ and
$I_{1-p}(\f^*) = I_p(\f)$).)

It is not hard to see (a proof is included below) that for every
$\epsilon >0$ there is a $C$ such that for every monotone $\f$ (on a
set of size $n$)

\beq{cpopt}
\mbox{$\exists ~p \in [n^{-\eps} p_c(\f),p_c(\f)]$ for which
$\f$ is $(C  \log (1/p) ,p)$-optimal.}
\enq
We don't know
whether $n^{-\eps}$ can be improved to $\Omega(1/\log n)$; in fact,
as far as we know even the following is possible.

\begin{conj}\label{C7}
For each $C>0$ there is an $\eps >0$ such that for any $\f$
$$
\mbox{$\exists ~p \in [\eps  p_c(\f)/\log n,p_c(\f)]$ for which
$\f$ is $(C  \log (1/p) ,p)$-optimal.}
$$
\end{conj}
Of course this
combined with Conjecture~\ref{C3}(c) would give
Conjecture~\ref{C1}.
(Considering that $\mu_p(\f)$ depends only on the sequence
$(a_0\dots a_n)$, where $a_k=a_k(\f):= |\{A\in\f:|A|=k\}|$,
and that possibilities for this sequence are completely
determined by
the Kruskal-Katona Theorem,
it's a little odd that Conjecture~\ref{C7}
should not be easy to settle
one way or the other.)

\mn
{\em Duality.}

\medskip
We would like to say---but don't know how---that even the $\log n$
``gap" of Conjecture~\ref{C1} is not typical and, more particularly,
implies some kind of simplicity in the dual.
Could it be, for instance, that if $p=p_c(\f) $ is bounded away from
0 and 1 (and the $\log n$ gap is tight),
then there is some $R\sub [n]$ of size $O(\log n)$ for which
$\mu_{1-p} (S \in \f ^*| S\supseteq R) > 1/2+\Omega(1)$
(equivalently, $\mu_p (S \in \f ^*| S\cap R=\0)
< 1/2-\Omega(1)$)?
And might this even extend to more general situations
if we replace $O(\log n)$
by $O(p^{-1}\log n)$?
%
A similar conclusion---existence of $R$ of size $O(\log n)$
for which either $\mu_p(\f|S\supseteq R)\approx 1$ or
$\mu_p(\f|S\cap R=\0)\approx 0$---is conjectured in
\cite[Conjecture 2.6]{Friedgut}, in case both $p_c(\f)$ and
$\mu_{p_c}(\f)$ are bounded away from 0 and 1,
and the isoperimetric inequality of \cite{KKL} is
tight to within a constant factor.

Slightly related thoughts about the
interplay of $\f$ and $\f^*$ may be found in \cite{FKW};
in particular, a variant of Conjecture 4.1 of that paper
would imply that for any {\em graph
property} $\f$ a weaker version of our Conjecture 1---in
which one only asks that $\langle \g\rangle$ contain most of
$\f$---holds for either $\f$ or $\f^*$.
In fact \cite{FKW}---which deals with
randomized decision tree complexity, especially of graph
properties---partly motivated our Conjecture 2;
but it is not clear how far this connection goes.
(For one thing, the aforementioned
Conjecture 4.1, and the complexity conjecture that motivated it,
are not correct for general properties.)

\mn {\em A little example}

\medskip
The following
example (the dual of the ``tribes" example of Ben-Or and
Linial \cite {BL}) may help to put some of the preceding discussion
in perspective.
Several of the above conjectures are based (to the extent they're
based on anything) on the idea that this simple construction is
about as bad as things get.

Suppose $k|n$, and let $X$ be a set of size $n$, and
$X_1\cup\cdots \cup X_m$ a partition of $X$, with $m=n/k$ and
$|X_i|=k$ for each $i$.
Let $\f$ consist of all sets meeting each $X_i$.
Then
$m(p)$ ($= \mu_p(\f)$) $ = (1-(1-p)^k)^{n/k} $ and
$
m'(p) = n (1-(1-p)^k)^{n/k-1}(1-p)^{k-1}.
$

If we take $k=\log n-\log\log n$, then $p_c(\f)$ is bounded away
from 0 and 1, while the conclusion of Conjecture 1
requires $q < O(1/\log n)$ (and we may then take $\g=\f$).
Here we achieve
$(O( \log (1/p)) ,p)$-optimality when $p=\theta (1/\log n)$,
while $(O(1),p)$-
and
$(1-o(1),p)$-optimality require, respectively,
$p = \log^{-(1+\Omega(1))} n$, and
$p = \log^{-\omega(1)} n$.

\mn {\em  Little proofs}

\mn {\em Proof of} (\ref{iso}).
We actually show
this for general (i.e. not necessarily monotone)
$\f$.  Note $I$ makes sense in this generality,
though its interpretation as a derivative is only valid
when $\f$ is monotone.

We proceed by induction on $n$.
Write $I'$, $\mu'$ $(=\mu_p')$ for influence and
measure in $\{0,1\}^{n-1}$.
Let $\A=\{x\in \f:x_n=0\}$, $\B=\{x\in \f:x_n=1\}$.
By induction,

$$
pI = p(1-p)I'(\A) +p^2I'(\B)+ p\mu'(\A\triangle \B)\geq
(1-p) \mu'(\A) \log_p\mu'(\A) + p\mu'(\B) \log_p\mu'(\B) +
p\mu'(\A\triangle \B),
$$
and we would like to say that the r.h.s. is at least

$$((1-p)\mu'(\A) + p\mu'(\B))\log_p((1-p)\mu'(\A) + p\mu'(\B)).$$
Set $\ga = \mu'(\A)$, $\gb =\mu'(\B)$.
Then assuming (w.l.o.g.) that $\gb\geq \ga$, it's enough to
verify the
numerical statement

$$(1-p)\ga\log_p\ga +p\gb\log_p\gb +p(\gb-\ga)-
((1-p)\ga +p\gb)\log_p ((1-p)\ga +p\gb) \geq 0.$$
But this holds with equality when $\gb=\ga$, and an easy calculation
shows that the derivative of the l.h.s. w.r.t. $\gb$ is nonnegative
(for {\em all} $\gb$, though we only need $\gb\geq\ga$).\qed

\mn
{\em Proof of} (\ref{cpopt}).
We prove this with $C \approx 1/(\eps \ln 2)$
(with the log in (\ref{cpopt}) interpreted as $\log_2$).
Set $p_c(\f)=p_c$.
Suppose $\f$ is not
$(C \log (1/p),p)$-optimal for any $p\in [n^{-\eps}p_c,p_c]$;
that is,
$$
pm'(p) > C m(p)  \log (1/m(p))
$$
for each such $p$.
We show this implies that $m(n^{-\eps}p_c)\leq 2^{-n}$,
a contradiction since
$p I =\sum\{\mu_p(x)|\{y\sim x:y\not\in\f\}|:x\in\f\} \leq n m(p)$
gives $(\log(1/p),p)$-optimality whenever $m(p)\leq 2^{-n}$.

Notice that
if $m(p) \le 2^{-r}$ and
$tm'(t) > C m(t)  \log (1/m(t))$
for $t\in [(1-1/(rC)) p, p]$,
then
$m(p( 1-1/(rC))) \le 2^{-(r+1)}$.
Thus our assumption gives
$m(p)\leq 2^{-n}$ for
$$p = p_c(1-1/C)(1-1/(2C))\cdots (1-1/((n-1)C)) \approx
n^{-\eps}p_c.$$\qed

\bn
{\bf Acknowledgment.}
We would like to thank Michel Talagrand for drawing our attention
to \cite{oldtal} and \cite{Talagrand}.

{\small

\begin{flushleft}
Department of Mathematics\\
Rutgers University\\
Piscataway NJ 08854 USA\\
jkahn@math.rutgers.edu
\end{flushleft}

\begin{flushleft}
Department of Mathematics\\
The Hebrew University\\
Jerusalem, Israel\\
kalai@math.huji.ac.il\\
and Departments of Computer Science and Mathematics, Yale
University
\end{flushleft}
}


\begin{thebibliography}{99}
\label{refs}


\bibitem {ANP} D. Achlioptas, A. Naor and Y. Peres
Rigorous location of phase transitions in hard optimization problems,
{\em Nature} {\bf 435} (2005), 759-764.


\bibitem{AS} N. Alon and J. Spencer,
{\em The Probabilistic Method}, Wiley, New York, 2000.


\bibitem {BL}
M. Ben-Or and N. Linial, Collective
coin flipping, in {\it Randomness and Computation} (S. Micali, ed.),
New York, Academic Press, pp. 91--115, 1990.

\bibitem{Boll81}
B. Bollob\'as, Random Graphs, pp. 257-274 in {\em Combinatorics},
London Math. Soc. Lecture Note Ser. 52, Cambridge Univ. Press,
Cambridge, 1981.

\bibitem{B}
B. Bollob\'as, The evolution of sparse graphs, pp. 35-57 in
{\em Graph Theory and Combinatorics}, B. Bollob\'as, ed.,
Academic Pr., 1984.


\bibitem{BRG}
B. Bollob\'as, {\em Random Graphs}, Academic Press, London, 1985.



\bibitem{BR}B. Bollob\'as and O. Riordan,
A short proof of the Harris-Kesten Theorem,
{\em Bull. London Math. Soc.}, to appear.



\bibitem{BT} B. Bollob\'as and A. Thomason, Threshold functions,
{\em Combinatorica} {\bf 7} (1987), 35-38.

\bibitem {Bou} J. Bourgain, On sharp thresholds of monotone properties,
Appendix to \cite {Fri}.



\bibitem{CFMR}
C. Cooper, A. Frieze, M. Molloy and B. Reed,
Perfect matchings in random $r$-regular, $s$-uniform hypergraphs,
{\em Combin. Probab. Comput.} {\bf 5} (1996), 1-14.

\bibitem{Erdos}
P. Erd\H{o}s,
On the combinatorial problems which I would most like to see solved,
{\em Combinatorica} {\bf 1} (1981), 25-42.


\bibitem{ER}
P. Erd\H{o}s and A. R\'enyi, On the evolution of random graphs,
{\em Publ. Math. Inst. Hungar. Acad. Sci.} {\bf 5} (1960), 17-61.

\bibitem{Fri} E. Friedgut,
Sharp thresholds of graph properties, and the $k$-sat problem.
{\em J. Amer. Math. Soc.} {\bf 12}  (1999), 1017-1054.

\bibitem{Friedgut}
E. Friedgut,
Influences in product spaces: KKL and BKKKL revisited,
{\em Combin. Probab. Comput.} {\bf 13} (2004), 17-29.


\bibitem{FKW} E. Friedgut, J. Kahn and A. Wigderson,
Computing graph properties by randomized subcube partitions,
{\em Randomization and Approximation Techniques
in Computer Science, 6th International Workshop} (2002), 105-113.

\bibitem{FJ}
A. Frieze and S. Janson,
Perfect matchings in random $s$-uniform hypergraphs,
{\em Random Structures \& Algorithms} {\bf 7} (1995), 41-57.

\bibitem{Grimmett}
G. Grimmett, {\em Percolation}, Second edition, Springer-Verlag,
Berlin, 1999.

\bibitem{JLR} S. Janson, T. \L uczak and A. Ruci\'nski,
{\em Random Graphs}, Wiley, New York, 2000.

\bibitem{Joh} A. Johansson,
Triangle factors of random graphs,
lecture at {\em Random Structures \& Algorithms},
Poznan, 2005.

\bibitem{KKL}
{ J. Kahn, G. Kalai and N. Linial }, The influence of variables
on Boolean functions, in
{\it Proc. 29-th Annual Symposium on Foundations of
Computer Science}, 68--80, 1988.

\bibitem {KaS}
G. Kalai and S. Safra,
Threshold Phenomena and Influence,  pp. 25-60 in
{\em Computational Complexity and Statistical Physics},
A.G. Percus, G. Istrate and C. Moore, eds.,
Oxford University Press, New York, 2006.



\bibitem{Kim} J.H. Kim,
Perfect matchings in random uniform hypergraphs,
to appear.

\bibitem{KS} J. Koml\'os and E. Szemer\'edi,
Limit distributions for the existence of Hamilton cycles in a random graph,
{\em Discrete Math.} {\bf 43} (1983), 55-63.

\bibitem{Korsh} A.D. Korshunov,
Solution of a problem of Erd\H{o}s and R\'enyi on Hamiltonian
cycles in non-oriented graphs, {\em Soviet Mat. Dokl.} {\bf 17}
(1976), 760-764.

\bibitem{Kriv} M. Krivelevich,
Triangle factors in random graphs,
{\em Combin. Probab. Comput.} {\bf 6} (1997), 337-347.


\bibitem {Ledoux}
M. Ledoux, The concentration of measure phenomenon,
Mathematical Surveys and Monographs, 89,
American Mathematical Society, Providence, RI, 2001.

\bibitem{Posa}  L. P\'osa, Hamiltonian circuits in random graphs,
{\em Disc. Math.} {\bf 14} (1976), 359-364.


\bibitem{Russo} L. Russo, On the critical percolation probabilities,
{\em Z. Wahrsch. Verw. Geb.} {\bf 56} (1981), 229-237.

\bibitem{SS}
J. Schmidt and E. Shamir,
A threshold for perfect matchings in random $d$-pure
hypergraphs, {\em Disc. Math.} {\bf 45} (1983), 287-295.


\bibitem{oldtal}
M. Talagrand,
Are all sets of positive measure essentially convex?, pp. 295--310 in
{\em Geometric Aspects of Functional Analysis (Israel, 1992--1994)}
(J. Lindenstrauss and V. Milman, eds.),
{\em Operator theory, advances and applications Vol. 77},
Bikh\"auser, Basel, 1995.


\bibitem{Talagrand}
M. Talagrand, Selector processes on classes of sets,
{\em Proba. Theor. Rel. Fields}, to appear.


\end{thebibliography}
\end{document}